\documentclass{amsart}[12pt]
\usepackage{amsmath, amsthm, amscd, amsfonts}

\makeatletter \oddsidemargin.9375in \evensidemargin \oddsidemargin
\theoremstyle{definition}

\theoremstyle{remark}

\numberwithin{equation}{section}

\begin{document}

\title{Differentiability Of Banach Spaces Via Constructible Sets }

\maketitle

\begin{center}
Hadi Haghshenas
\end{center}

\begin{center}
Department of Mathematics, Birjand University, Iran\\
E-mail Address: h$_{-}$haghshenas60@yahoo.com
\end{center}

\begin{center}
\textbf{Abstract }
\end{center}
The main goal of this paper is to prove that any Banach space $X$,
that every dual ball in $X^{**}$ is $weak^{*}-separable$, or every
$weak^{*}-closed$ convex subset in $X^{**}$ is
$weak^{*}-separable$, or every norm-closed convex set in $X^{*}$
is constructible, admits an equivalent Frechet differentiable
norm.

\textbf{Key Words:}\hspace{.25cm} Frechet differentiability,
locally uniformly convex norm, biorthogonal system, countable
intersections, constructible sets.

\section{\textbf{INTRODUCTION}}
We consider a real Banach space $X$, and shall call a closed
convex set $ C\subset X $ constructible if, it is the countable
intersection of closed half-spaces in $X$. Recall that an open \{
resp. A closed\} Half-space of $X$ is a set of the form
$H=f^{-1}(-\infty, \alpha)\{ resp.  H=f^{-1}(-\infty, \alpha ] \}$
with $f \in X^{*}\backslash \{ 0 \}$ and $\alpha \in \Bbb{R} $. We
will begin with some basic facts concerning constructible sets,
many of which we will use frequently
without explicit mention.\\
\textbf{Theorem 1. }[2] (a) Any translate of a constructible set
is constructible.\\(b) If $C$ is constructible, then $\lambda C$
is constructible for all $ \lambda\neq 0 $. \\ (c) A nonempty
intersection of countably many constructible sets is
constructible.\\(d) If $0 \in C$ and $C$ is constructible, then we
can represent, $$ C= \displaystyle{\ \bigcap_{n = 1}^{\infty}
f_{n}^{-1}(\ -\infty , 1 ] }$$\\ We give a simple criterion which
we will use to build closed convex sets that are not
constructible.\\\textbf{Theorem  2. }Let $C$ be a closed convex
subset of a Banach space $X$. Suppose there is an uncountable
family $\{x_{i}\}_{i \in I}\subset X$ such that $x_{i}$ not in $C$
for all $i \in I$, but $(x_{i} + x_{j}) /2 \in
C $ for all $i\neq j$. Then $C$ is not constructible.\\
\textit{proof .   }By translation we may assume $0 \in C$. Now let
us suppose $ C= \displaystyle{\ \bigcap_{n = 1}^{\infty}
f_{n}^{-1}(\ -\infty , 1 ] }$. For each $i$, we choose $n_{i}$such
that $f_{n_{i}}(x_{i}) > 1$. Now $(x_{i} + x_{j}) /2 \in C $ and
so $f_{n_{i}}(x_{j}) < 1$ for all $i\neq j$. Therefore, if $i\neq
j$ , one has $f_{n_{i}}(x_{i}) > 1$ and $f_{n_{j}}(x_{i}) < 1$.
This shows $n_{i} \neq n_{j}$ for $i\neq j$ , which is a
contradiction.\section{SEPARABILITY}To approach differentiability
of a Banach space $X$, we will obtain conditions that are
equivalent by separability of $X^{*}$. Then we will show that
every Banach space with a separable dual, has an equivalent
Frechet differentiable norm.\\\textbf{Theorem 3. }[2] Let $X$ be a
Banach space, then the following are equivalent,\\(a) There is an
uncountable family$\{x_{\alpha}\}\subset X$ such that $x_{\alpha}$
not in $\overline{conv}(\{x_{\beta}:\beta\neq\alpha\})$ for all
$\alpha$.\\(b) There is a bounded closed convex subset of $X$ that
is not constructible.\\(c) There is a closed convex subset in $X$
that is not constructible.\\(d) There is a $weak^{*}-closed$
convex subset of $X^{*}$that is not $weak^{*}-separable$.\\(e)
There is a ball of an equivalent dual norm in $X^{*}$that is not
$weak^{*}-separable$.\\(f) There is an equivalent norm on $X$
whose unit ball is not constructible.\\(g) There is an uncountable
biorthogonal system in $X^{*}$.\\\\There are several other
conditions equivalent to
those listed in theorem 3, that can be found in [4].\\
\textbf{Remark  4. }Suppose $X$ is a nonseparable dual Banach
space . Then $X$ admits an uncountable biorthogonal
system.\\\textit{proof .   }We write $X=Z^{*}$ for some Banach
space $Z$. If $Z^{*}$ has the Radon-Nikodym property, then $Z^{*}$
admits an uncountable biorthogonal system according to [3].
Otherwise, there is a separable subspace $Y \subset Z $ such that
$Y^{*}$ is not separable; by [5, corollary 3], $Y^{*}$ has an
uncountable biorthogonal system. Thus, $X^{*} / Y^{\bot}$ has an
uncountable biorthogonal system, which can be easily pulled back
to $X^{*}$.\\\textbf{Corollary  5. }A dual space $X^{*}$ is
separable, if and only if, every norm - closed convex set in
$X^{*}$ is constructible.\\\textit{proof .   }If $X^{*}$ is
separable, then every closed convex subset of $X^{*}$ can be
written as a countable intersection of half-spaces by [1].
Conversely, if $X^{*}$ is not separable, there is an uncountable
biorthogonal system in $X^{*}$. (as cited in remark 4 ). Thus
theorem 3 ensures that there is a closed convex set in $X^{*}$
that is not constructible.\\\\As a further observation, we provide
a characterization of separable Asplund spaces in terms of
$weak^{*}-separability$ of $weak^{*}-closed$ convex sets in the
second dual.\\\textbf{Corollary  6. }For a separable Banach space
$X$, the following are equivalent,\\(a) $X^{*}$ is separable.\\(b)
Every dual ball in $X^{**}$ is $weak^{*}-separable$.\\(c) Every
$weak^{*}-closed $ convex subset in $X^{**}$ is
$weak^{*}-separable$.\\ \textit{proof .   } $(a) \Rightarrow (c)$
follows from the separability of $X^{*}$, and $(c) \Rightarrow
(b)$ is trivial. To prove $(b) \Rightarrow (a)$ we suppose $X^{*}$
is not separable. According to corollary 5 and theorem 3, there is
an equivalent dual ball in $X^{**}$ that is not
$weak^{*}-separable$.\\\section{\textbf{DIFFERENTIABILITY}}The
norm function $\|.\|$ is Gateaux differentiable at $0\neq x$ if,
there is some $x^{*}\in X^{*}$ such that the quantity
$\displaystyle{\lim_{t\rightarrow 0}\frac{\|x+ty\|-\|x\|}{t}}$
exists for each $y\in X$ and is equal to $\langle x^{*},y
\rangle$. If the limit exists uniformly for each $y$ in unit
sphere $S(X)$, then the norm function is said to be Fr\'{e}chet
differentiable at $x$.
\\\textbf{Theorem 7. }The norm function is Fr\'{e}chet
differentiable at $0\neq x\in X$, if and only if,
$$\displaystyle{\lim_{t\rightarrow 0}\frac{\|x+ty\|+\|x-ty\|-2\|x\|}{t}=0,}$$
uniformly for each $y\in S(X)$.\\\textbf{Theorem 8. }If $X^{*}$ is
separable, then $X$ has an equivalent norm wich is Frechet
differentiable at every $x\neq 0$.\\\textit{proof . }Let
$\{x_{i}\}_{i=1}^{\infty}$ be a dense sequence in the unit ball of
$X$. the linear operator $T:X^{*} \rightarrow l_{2}$ with
$Tu^{*}=(\langle u^{*},x_{1} \rangle , ... , \langle
u^{*},x_{n}\rangle/n,...)$ is bounded and injective from $X^{*}$
into $l_{2}$. Fix an increasing sequence
$\{F_{n}\}_{n=1}^{\infty}$ of finite-dimensional subspaces of
$X^{*}$ whose union is norm dense in $X^{*}$. Define a nonlinear
homogeneous map $S:X^{*} \rightarrow l_{2}$ by $$S(u^{*})=(
d(u^{*},F_{1}), ... ,d(u^{*},F_{n})/n )$$where the distances are
with respect to the original norm in $X^{*}$. Then the components
of $S$ are $w^{*}-lower$ $ semi-continuous $ (because the
$F_{n}^{,}$s are finite-dimensional) and $\|S(u^{*}+ v^{*})\|_{2}
\leq \|S (u^{*})\|_{2}+\|S( v^{*})\|_{2}$. Thus the norm $\||.|\|$
on $X^{*}$ defined by $$\||u^{*}|\|^{2}=\|u^{*}\|^{2}+\|T
u^{*}\|_{2}^{2}+\|S( u^{*})\|_{2}^{2}$$ is an equivalent norm with
a $w^{*}-closed$ unit ball ; i.e., it is dual to an equivalent
norm on $X$ ( also denoted by $\||.|\|$ ). We first show that
$(X^{*}, \||.|\|)$ is locally uniformly convex; i.e., if $u^{*},
u_{k}^{*}\in X^{*}$satisfy $\||u^{*}|\|=\||u_{k}^{*}|\|=1$ and
$\displaystyle{\lim_{k\rightarrow \infty}\||u^{*}+u_{k}^{*}|\|=2}$
, then $\||u_{k}^{*}-u^{*}|\| \longrightarrow 0$.\\Indeed, each of
the three expressions whose sequares appear in definition of
$\||.|\|$ satisfy the triangle inequality and all three are
$w^{*}-lower$ $ semi-continuous $. It follows that
$\|u_{k}^{*}\|\longrightarrow \|u^{*}\|$ and similarly for the
other two terms. A similar argument shows that
$\displaystyle{\lim_{k\rightarrow
\infty}\|T(u^{*}+u_{k}^{*})\|_{2}=2 \|Tu^{*}\|_{2}}$ and thus
$\|T(u_{k}^{*}-u^{*})\|_{2} \longrightarrow 0$. Since
$\{x_{i}\}$is dense in the unit ball of $X$, it follows that
$u_{k}^{*}$ tends to $u^{*}$ in the $w^{*}-topology$. By
considering $S$ we obtain similarly
$\displaystyle{\lim_{k\rightarrow
\infty}d(u_{k}^{*},F_{n})=d(u^{*},F_{n})}$ for each fixed $n$. Let
$\varepsilon > 0$, and pick $n$ such that $d(u^{*},F_{n}) <
\varepsilon $. Then $d(u_{k}^{*},F_{n})< \varepsilon $ for large
enough $k$, so find $z_{k}^{*} \in F_{n}$ such that
$\|u_{k}^{*}-z_{k}^{*}\|=d(u_{k}^{*},F_{n}) < \varepsilon $.
Assume ( by passing to a subsequence ) that the $z_{k}^{*}$$^{,}$s
converge to a point $z^{*}$. Then $\|u_{k}^{*}-z^{*}\|\leq
2\varepsilon $ for large $k$, and $\|u^{*}-z^{*}\|\leq
2\varepsilon $ because $u_{k}^{*}$ is $w^{*}-convergent$to
$u^{*}$.Thus $lim sup \|u_{k}^{*}-u^{*}\|\leq 4\varepsilon $, and
since $\varepsilon$ was arbitrary
$\|u_{k}^{*}-u^{*}\|\longrightarrow 0 $.\\Now, we show that
whenever a dual norm $\||.|\|$ on a dual space $X^{*}$ is locally
uniformly convex, then the norm on $X$ is Frechet differentiable
at every $x \neq 0$. Indeed, assume that $\||x|\|=1$, and choose
$u^{*} \in X^{*}$with $\||u^{*}|\|=\langle u^{*}, x \rangle =1 $.
For every $y \in X $ choose norm-one functionals $v_{y}^{*},
w_{y}^{*} \in X^{*}$ such that $\langle v_{y}^{*},
x+y\rangle=\||x+y|\|$ and $\langle
w_{y}^{*},x-y\rangle=\||x-y|\|$. Then $$\||v_{y}^{*}+u^{*}|\| \geq
\langle v_{y}^{*} + u^{*}, x \rangle = 1 + \||x+y|\|- \langle
v_{y}^{*}, y\rangle \longrightarrow 2 $$ as $y \longrightarrow 0$.
By local uniform convexity $\||v_{y}^{*}-u^{*}|\|\longrightarrow 0
$ as $y \longrightarrow 0$, and
$\||w_{y}^{*}-u^{*}|\|\longrightarrow 0 $ in a similar way. Thus
$\||v_{y}^{*}-w_{y}^{*}|\|\longrightarrow 0 $, and hence
\begin{eqnarray*}
\||x+y|\|+ \||x-y|\|-2 & = & \langle v_{y}^{*} , x \rangle +
\langle w_{y}^{*} , x \rangle + \langle v_{y}^{*}- w_{y}^{*}, y\rangle-2
\\& \leq & \||v_{y}^{*}-w_{y}^{*}|\| \||y|\| = o(\||y|\|).
\end{eqnarray*}

\end{document}